# Two Iterative Algorithms for Solving Systems of Simultaneous Linear Algebraic Equations with Real Matrices of Coefficients

A. S. Kondratiev (1 and 2) and N. P. Polishchuk (2)

(1) Moscow Power Engineering Institute,
(2) Altair Naval Research Institute of Radio Electronics

**Abstract:** The paper describes two iterative algorithms for solving general systems of $M$ simultaneous linear algebraic equations (SLAE) with real matrices of coefficients. The system can be determined, underdetermined, and overdetermined. Linearly dependent equations are also allowed. Both algorithms use the method of Lagrange multipliers to transform the original SLAE into a positively determined function $F$ of real original variables and Lagrange multipliers $\lambda_m$. Function $F$ is differentiated with respect to variables $x_i$ and the obtained relationships are used to express $F$ in terms of Lagrange multipliers $\lambda_m$. The obtained function is minimized with respect to variables $\lambda_m$ with the help of one of two the following minimization techniques: (1) relaxation method or (2) conjugate gradient method by Fletcher and Reeves. Numerical examples are given.

## INTRODUCTION

Most of well-known methods for solving systems of simultaneous (consistent) linear algebraic equations that can be found in standard software packages either fail to solve systems with degenerate matrices of coefficients (for the systems containing linearly dependent equations) and the systems in which the number of equations does not coincide with the number of unknowns (in the case of underdetermined or overdetermined systems) or require tedious preliminary manipulations. Therefore, the development of a universal algorithm that is suitable for solving both underdetermined and overdetermined systems and does not require a preliminary analysis of the system type seems to be useful.

## 1. FORMULATION OF THE PROBLEM

Let us we have a system of simultaneous linear equations

$$[A]x> = b>, \qquad (1)$$

where $x>$ is the column vector of $N$ real unknown variables $x_n$, $[A]$ is the $M$x$N$ real matrix of coefficients, $M$ is the number of equations, $N$ is the number of unknowns, and $b>$ is the real column vector of length $M$.

It is required to find the least-mean-square solution to system (1). If $M \leq N$, system (1) can be solved exactly; otherwise, the solution is approximate.

## 2. SOLUTION TECHNIQUE

This problem is reduced to the constrained minimization problem for the quadratic functional

$$Q = Q(x_1,\ldots, x_N) = \sum_{n=1}^{N} x_n^2 \qquad (2)$$

with constraints (1).

Using the method of Lagrange multipliers [1], the constrained minimization problem (1), (2) transforms to an unconstrained minimization problem for the following functional:

$$\Phi = \Phi(x_1,\ldots, x_N, \lambda_1,\ldots, \lambda_M) = Q + <\lambda [A] x> - <\lambda\, b>, \quad (3)$$

where $<\lambda$ is the column vector of length $M$ composed of Lagrange multipliers $\lambda_m$ ($m = 1,\ldots, M$).

Differentiating functional (3) with respect to $x_n$ and equating the obtained expressions to zero, we find relationships between desired unknowns $x_n$ and Lagrange multipliers $\lambda_1,\ldots, \lambda_M$:

$$x_n = -0.5 <\lambda a_n>, \quad n = 1,\ldots, N, \quad (4)$$

where $a_n>$ is the $n$th column of matrix $[A]$.

Substituting formulas (4) into (3), we can express functional $\Phi$ as a quadratic form depending only on Lagrange multipliers $\lambda_1,\ldots, \lambda_M$:

$$\Phi(\lambda_1,\ldots, \lambda_M) = -0.25 <\lambda [A][A]_t \lambda> - <\lambda\, b>, \quad (5)$$

where subscript $t$ denotes transposition.

Therefore, the problem reduces to minimization of function $\Phi(\lambda_1,\ldots, \lambda_M)$ treated as a function of $M$ variables $\lambda_1,\ldots, \lambda_M$.

The minimum of functional $\Phi$ is sought using either of two the following iterative procedures: (1) the relaxation (coordinate-by-coordinate) minimization method [2] or (2) the conjugate gradient method by Fletcher and Reeves [3].

This approach allows us, first, to decrease the number of unknowns to the number of the equations in system (2) and, second, to avoid application of direct methods for solving the SLAE corresponding to functional (5), which can be important for large equation systems. It also does not involve an explicit analysis of the type of the solved system.

In the case of minimization with the use of the relaxation method, the iterative procedure contains the following steps:
1. Set the initial point as $<\lambda_k = 0$, where $k = 0$ ($k$ is the step of the iterative procedure).
2. At the $k$th step ($k = 1,\ldots, M$), we fix all variables $\lambda_1,\ldots, \lambda_{k-1}, \lambda_{k+1},\ldots, \lambda_M$ and seek the minimum of functional $\Phi$ treated as a function of only one variable $\lambda_k$:

$$\lambda_k = (2b_k + p_{kk}) / w_{kk}, \quad (6)$$

where

$$p_{kk} = <\lambda w_k> - \lambda_k^{old} w_{kk},$$
$$w_{lm} = \sum_{n=1}^{N} a_{lm} a_{mn}, \quad (7)$$

$w_k>$ is the $k$th column of matrix $[w]$ composed of elements $w_{lm}$, and $\lambda_k^{old}$ is the value of the $k$th variable found at the preceding step.
Note that $[w] = [A][A]_t$.
3. If, upon finishing the passage over all variables, the sum of absolute values of changes in all $M$ variables is found to be less that some prescribed value, we assume that the minimum of functional (5) was reached. Otherwise, we return to step 2.

In the case of minimization with the use of the conjugate gradient method, the iterative procedure contains the following steps:
1. Set the initial point as $\lambda_i> = 0$, where $i = 1$ ($i$ is the iteration number).
2. Calculate the gradient of the minimized functional

$$r_i >= -0{,}5 < \lambda_i[w] - b >. \tag{8}$$

3. Calculate the absolute value of the gradient and compare it with some prescribed accuracy value. If the absolute value of the gradient is less than the prescribed small value, the minimization procedure terminates; otherwise, it passes to to Step 4.
4. Determine the search direction for the minimization of functional (5).
If $i = 1$ (the first iteration), the search direction is specified by the column vector $d_i>$ expressed as

$$d_i >= -r_i >. \tag{9a}$$

At subsequent $M$ iterations ($i = 2,\ldots, M+1$),

$$d_i >= -r_i > + \beta_{i-1} \cdot d_{i-1} >, \tag{9b}$$

where coefficient $\beta_{i-1}$ is determined from the formula

$$\beta_{i-1} = \frac{< r_i r_i >}{< r_{i-1} r_{i-1} >}, \tag{10}$$

5. Determine the optimum point from three vales of the minimized quadratic functional $\Phi$. For this purpose, calculate values of functional $\Phi$ at three points: the current point ($\Phi_2$) and two points ($\Phi_1$ and $\Phi_3$) equispaced from it by distances $\pm\delta$ along the direction determined by formulas (9):

$$\lambda_i >= \lambda_{i-1} > + \frac{(\Phi_1 - \Phi_3)\delta}{2(\Phi_1 - 2\Phi_2 + \Phi_3)} d_i > \tag{11}$$

6. Compare the sum of absolute values of all variables $\lambda_{k,i}$ ($k$ is the iteration number) with some prescribed accuracy value and either terminate iterations (if the calculated sum is less than the prescribed value) or return to Step 2.

Theoretically [3], the conjugate gradient method can find the minimum of a quadratic functional after $M+1$ iterations. However, a finite computational accuracy requires multiple passages over all variables to ensure the specified accuracy. This fact is illustrated by the examples of calculations presented below.

In both SLAE solution algorithms, desired unknown variables $x_n$ are determined from the values of Lagrange multipliers $\lambda_k$ found from the iterative procedure with the help of formula (4).

### 3. NUMERICAL EXAMPLES

1. Determined SLAE with $N = M = 3$; $\det[A] \neq 0$.

a)
$$x_1 + x_2 = 2.0$$
$$x_2 + x_3 = 2.0$$
$$x_1 + x_3 = 2.0$$

Both algorithms yielded the same result obtained after one passage over $M$ variables: $x_1 = 1.0$; $x_2 = 1.0$; $x_3 = 1.0$

b)
$$33x_1 + 16x_2 + 72x_3 = 129.0$$
$$-24x_1 - 10x_2 - 57x_3 = -96.0$$
$$18x_1 - 11x_2 + 7x_3 = 8.5$$

The exact solution is $x_1 = 1.0$; $x_2 = 1.5$; $x_3 = 1.0$

The first algorithm gives the solution $x_1 = 1.0000000$; $x_2 = 1.5000000$; $x_3 = 1.000000$ which was obtained after 7098 passages over $M$ variables.

The second algorithm gives the solution $x_1 = 0.9999183$; $x_2 = 1.499904$; $x_3 = 1.000052$. This solution was obtained after 2657 passages over $M$ variables.

2. Degenerate SLAE containing two linearly dependent equations ($N = M = 3$; $\det[A]$
$$x_1 + x_2 + x_3 = 1.0$$
$$x_1 + x_2 + x_3 = 1.0$$
$$x_1 - x_2 = 0.$$

Both algorithms came to the same result (obtained after a single passage):
$x_1 = 0.33333$; $x_2 = 0.33333$ ; $x_3 = 0.33333$.

3. Overdetermined SLAE with an ill-conditioned matrix of coefficients ($M > N$).

$$x_1 + 2x_2 + 4x_3 = 4.999$$
$$x_1 + 4x_2 + 16x_3 = 9.001$$
$$x_1 + 6x_2 + 36x_3 = 12.999$$
$$x_1 + 8x_2 + 64x_3 = 17.001$$

The solution given in [4] is
$x_1 = 0.998997$; $x_2 = 2.000200$; $x_3 = -2.8825692 \cdot 10^{-8}$.
The first algorithm comes to the following solution:
$x_1 = 0.996995$; $x_2 = 2.001182$; $x_3 = -8.5153115 \cdot 10^{-5}$.
The solution obtained with the second algorithm is
$x_1 = 0.999000$; $x_2 = 2.000204$; $x_3 = +2.3505518 \cdot 10^{-5}$.

The results presented above were obtained after 250000 passages over all variables for a specified accuracy for termination of the iterative procedure of $10^{-10}$.

**CONCLUSION**

The results presented above confirm the efficiency af the described algorithms that can be used for solving systems of simultaneous linear equations without any preliminary analysis of the matrix of coefficients. Numerical experiments have shown that the first algorithm can even find approximate solutions to inconsistent systems of linear equations.

The approach described can be straightforwardly extended to the solution of systems with complex matrices of coefficients.

**REFERENCES**

1. G. Korn and T. Korn, *Mathematical Handbook for Scientists and Engineers* (McGraw-Hill, New York, 1968).
2. M. Aoki, *Introduction to Optimization Techniques: Fundamentals and Applications of Nonlinear Programming*. (Macmillan, New York, 1971).
3. R. Fletcher, C.M. Reeves, Function minimization by conjugate gradients, *Comp. J.* **7**, 149 (1964).
4. *IMSL Math/Library: Fortran Subroutines for Mathematical Applications* (Visual Numeric, Houston (Texas), 1997).